\documentclass[a4paper,12pt]{article}

\usepackage{graphicx}      
\usepackage{epsfig} 
\usepackage{amsmath} 
\usepackage{amssymb}  

\newcommand{\RR}{\mathbb{R}}
\newcommand{\NN}{\mathbb{N}}

\newcommand{\D}{\mathcal{D}}

\newcommand{\nocode}[1]{}

  \newtheorem{prop}{Proposition}

  \newtheorem{method}{Method}

\title{\LARGE\bf
Moment and SDP relaxation techniques for smooth approximations of problems involving nonlinear differential equations} 

\begin{document}

\author{Martin Mevissen$^1$, Jean B. Lasserre$^2$, Didier Henrion$^{2,3}$}

\footnotetext[1]{Department of Mathematical and Computing Sciences,
Tokyo Institute of Technology, Ookayama 2-12-1, Meguro-ku, Tokyo 152-8552, Japan.}

\footnotetext[2]{CNRS; LAAS; 7 avenue du colonel Roche, F-31077 Toulouse, France;
Universit\'e de Toulouse; UPS, INSA, INP, ISAE; LAAS; F-31077 Toulouse, France.}

\footnotetext[3]{Faculty of Electrical Engineering,
Czech Technical University in Prague, Technick\'a 4, CZ-16626 Prague, Czech Repu
blic.}

\maketitle

\begin{abstract}                
Combining recent moment and sparse semidefinite programming (SDP) relaxation techniques, we propose an approach to find smooth approximations for solutions of problems involving nonlinear differential equations. Given a system of nonlinear differential equations, we apply a technique based on finite differences and sparse SDP relaxations for polynomial optimization problems (POP) to obtain a discrete approximation of its solution. In a second step we apply maximum entropy estimation (using moments of a Borel measure associated with the discrete solution) to obtain a smooth closed-form approximation.
The approach is illustrated on a variety of linear and nonlinear ordinary differential equations (ODE), partial differential equations (PDE) and optimal control problems (OCP), and preliminary numerical results are reported.
\end{abstract}

\begin{center}{\bf Keywords:}
semidefinite programming, nonlinear optimal control, maximum entropy estimation, approximative methods, nonlinear optimization, moment methods.
\end{center}


\section{Introduction}
\label{secIntro}
Problems involving nonlinear differential equations arise in a variety of models for real world problems. Even finding {\it approximate} solutions for nonlinear differential equations 
remains a challenge. In the previous work \cite{mknt,myt} we have established a technique based on sparse semidefinite programming (SDP) relaxations to construct discrete approximations for solutions of systems of nonlinear differential equations. In this paper we present a novel approach to obtain smooth approximations for solutions of differential equations and of problems involving differential equations such as optimal control problems. Namely, an approximate solution is obtained by applying the maximum entropy estimation of \cite{borwein,lassCDC} with a finite number of moments of a Borel measure
associated with a discrete approximation of the solution of the differential equation. For
{\it linear} differential equations, an SDP relaxation based method was proposed in \cite{bertsimasCaramanis} to generate contracting sequences of lower and upper bounds for the moments.

Our contribution in this paper is primarily concerned with {\it nonlinear} differential equations and {\it nonlinear} optimal control problems.
In a first step, we take advantage of the discrete approximations provided by the SDP relaxation method in \cite{mknt} to compute a finite set of moments for an appropriately defined measure with discrete support.
Next, an approximation for the solution of the differential equation is obtained in closed form by maximum entropy estimation, using the moments of the discrete measure. To the best of our knowledge, it seems to be the first attempt to apply maximum entropy estimation to obtain smooth approximations for solutions of linear and nonlinear differential equations. Finally, if maximum entropy estimation does guarantee some weak convergence of the estimate to the true solution as the number of moments increases,
it does not guarantee {\it pointwise} convergence on the entire domain of the differential equation. However, as our preliminary results show for different linear and nonlinear differential equations and optimal control problems, accurate pointwise approximations can be achieved on certain regions of the domain.

The structure of the paper is as follows. In Section \ref{secSDPR} we briefly recall the sparse SDP relaxations method for solving nonlinear differential equations numerically, which is the basis for our technique. In Section \ref{secMaxEntropy} we introduce the method of maximum entropy estimation. Our technique to compute smooth approximations for the solutions of
ordinary differential equations (ODE), partial differential equations (PDE) and optimal control problems (OCP), which combines the methods from Section \ref{secSDPR} and \ref{secMaxEntropy}, is presented in Section \ref{secSmooth}. Finally, preliminary numerical results for this method are reported in Section \ref{secNumExp}.

\section{Sparse SDP relaxations for solving differential equations}
\label{secSDPR}

In this section we recall the approach to compute discrete approximations for solutions to systems of differential equations with polynomial data presented in \cite{mknt}.

\subsection{Transforming a differential equation into a POP}
In this paper we consider bidimensional differential equation problems of the following type
\begin{equation}
\begin{array}{rl}
L(u(x,y)) + G(u(x,y)) = f(x,y) & \forall (x,y)\in \Omega,\\
H(u(x,y)) = g(x,y) & \forall (x,y) \in \partial\Omega,\\
\text{lbd} \leq u(x,y) \leq \text{ubd} & \forall (x,y)\in \Omega,
\end{array}
\label{PDEproblem}
\end{equation}
with $\Omega:=[x_{\min},x_{\max}]\times [y_{\min},y_{\max}] \subset \mathbb{R}^2$, $f:\, \Omega\rightarrow\RR^m$, $g:\, \partial\Omega\rightarrow\RR^m$ smooth functions, $L:\,\D\rightarrow\D$ and $H:\,\D\rightarrow\D $ differential operators, $G: \, \D\rightarrow\D$ operator, $L$, $H$ and $G$ are polynomial in the function $u\in\D$ and its derivatives for some function space $\D\subseteq L_1(\Omega)$, with $L_1(\Omega)$ the Banach space of integrable functions on $\Omega$, and $\text{lbd}<\text{ubd}$ bounds for the function values of $u$. Even if lower and upper bounds are not given by the problem directly, we need to impose them, since they are crucial for numerical stability of our method. If little is known about the solutions of (\ref{PDEproblem}), loose lower and upper bounds are to be imposed.

In a first step of the sparse semidefinite programming relaxation (SDPR) method for solving differential equations,  we transform an ODE or a PDE of form (\ref{PDEproblem}) into a polynomial optimization problem (POP). We then discretize the rectangular domain $\Omega$ by $N$ (one-dimensional case) or $N_x N_y$ (two-dimensional case) grid points, approximate the derivates by standard finite differences, for instance $\frac{\partial^2 u(x_i,y_j)}{\partial x^2} \approx \frac{u_{i+1,j} - 2u_{i,j}+u_{i-1,j}}{{\Delta x}^2}$, where $u_{i,j}:=u(x_i,y_j)$, and denote the discretized differential operators at some interior grid point $(x_i,y_j)$ and boundary grid point $(x_k,y_l)$ as $L_{i,j}$ and $H_{k,l}$, respectively. Finally, we take the discretization of the problem (\ref{PDEproblem}) as a system of constraints, choose some objective function $F$ which needs to be polynomial in $u$, i.e. $F(u)\in\RR[u]$, and obtain the following optimization problem
\begin{equation}
\begin{array}{l@{\;}l}
\min & F(u)\\
\text{s.t. } & L_{i,j}(u) + G(u_{i,j}) = f(x_i,y_j), \\ 
& \text{lbd}_{i,j} \leq u_{i,j} \leq \text{ubd}_{i,j} \:\:  \forall  (i,j) \in\{1..N_x\} \times \{ 1..N_y\},\\
& H_{k,l}(u) = g(x_k,y_l)  \:\:\forall (k,l) \in \{1,N_x\}\times \{1..N_y\} \\
& \quad \cup \{ 1..N_x \}\times \{ 1,N_y \} \\
\end{array}
\label{discretizedPDEPOP}
\end{equation}
which, since all functions are polynomial in the variable $u=(u_{1,1},\ldots,u_{N_x,N_y})$, is a POP. In the case of an OCP the objective function $F$ is given by the discretization of the optimal value function. In the case of a differential equation with several solutions we can pick a specific solution by choosing an appropriate objective function. A particular choice for the objective function may be motivated by the underlying problem, such as energy functionals in physics.

\subsection{Sparse SDP relaxations}
Problem (\ref{discretizedPDEPOP}) is of the form:
\begin{equation}
\begin{array}{llll}
\text{POP} & \min & F(u) \\
& \text{s.t.} & \tilde{g}_j(u) \geq 0 & \forall j\in\{1,\ldots,\tilde{k}\},\\
&& \tilde{h}_i(u) = 0 & \forall i\in \{ 1,\ldots,\tilde{l} \},\\
&& \text{lbd}_s \leq u_s \leq \text{ubd}_s & \forall s\in\{1,\ldots,\tilde{n}\},
\end{array}
\label{genPOP}
\end{equation}
where $F(u)=\sum_{\alpha\in\NN^{\tilde{n}}}F_{\alpha}u^{\alpha}$, $\tilde{g}_j(u) = \sum_{\alpha\in\NN^{\tilde{n}}}{g_j}_{\alpha}u^{\alpha}$ and $\tilde{h}_i(u) = \sum_{\alpha\in\NN^{\tilde{n}}}{h_i}_{\alpha}u^{\alpha}$. 
As mentioned above, the dimension of this POP is $\tilde{n}=N_x\, N_y$ (or $\tilde{n}=N$ in the one-dimensional case). But due to the structure of the finite difference discretization only a small number of the $\tilde{n}$ components of $u$ occurs in each constraint or in each monomial of the objective function of (\ref{discretizedPDEPOP}). Thus, the POPs derived from nonlinear differential equations are sparse. 

A systematic way of characterizing the structured sparsity of a POP has been introduced in \cite{wkkm}, and this structured sparsity is exploited by the {\it sparse SDP relaxations} constructed in \cite{wkkm}. For the POP (\ref{genPOP}) one obtains the sparse SDP relaxations:
\begin{equation}
\begin{array}{llll}
\text{SDP}_w & \min & \sum_{\mid\alpha\mid\leq 2w} F_{\alpha}y_{\alpha}\\
& \text{s.t.} & M_{w-w_j}(\tilde{g}_j\, y, I_{t(j)}) \succcurlyeq 0 & \forall\, j\in \{1,\ldots,\tilde{k}\}\\
&& M_{w-\tilde{w}_i}(\tilde{h}_i\, y, I_{\tilde{t}(i)}) = 0 & \forall\, i\in \{1,\ldots,\tilde{l}\}\\
&& M_w(y, I_t) \succcurlyeq 0 & \forall \, t\in \{1,\ldots,d\}\\
&& \text{lbd}_s \leq y_s \leq \text{ubd}_s & \forall \, s\in\{1,\ldots,\tilde{n}\},  
\end{array}
\label{sparseSDP}
\end{equation}
where $w_j:=\lceil\frac{\deg \tilde{g}_j}{2}\rceil$, $\tilde{w}_i:=\lceil\frac{\deg \tilde{h}_i}{2}\rceil$, $\{I_t\}_{t=1}^d$ is the set of subsets of $\{1,\ldots,\tilde{n}\}$ derived from the correlative sparsity pattern matrix of the POP \cite{wkkm}, $M_w(\cdot,\cdot)$ and $M_{w-w_j}(\cdot,\cdot)$ are the partial moments and localizing matrices of \cite{lassSparse}, and $w\geq w_{\min}:=\max_{i,j}(w_j,\tilde{w}_i)$ is the order of the SDP relaxation. 
As $w$ increases, solving $\text{SDP}_w$ generates a nondecreasing sequence of lower bounds for $\min(\text{POP})$,
namely:
\begin{equation*}
\min \left(\text{SDP}_{w_{\min}}\right) \leq \min \left(\text{SDP}_{w_{\min}+1}\right) \leq \ldots \leq \min \left(\text{POP}\right).
\end{equation*}
Moreover, if (\ref{genPOP}) has a unique minimizer $x^{\star}$, the vector
$(y_{w,1}^{\star},\ldots,y_{w,{\tilde{n}}}^{\star})$ obtained from an optimal solution $y_w^{\star}$ of (\ref{sparseSDP}) is an approximation of $x^{\star}$, and under certain compactness conditions, see \cite{lassSparse}, as $w\to\infty$,
\begin{equation}
\min\left( \text{SDP}_w \right) \rightarrow \min\left(\text{POP}\right) \text{ and } y_w^{\star} \rightarrow x^{\star}.
\label{conv}
\end{equation}
In implementations of the sparse SDP relaxation such as SparsePOP, small linear perturbation terms are added to the objective function $F$, in order to ensure that (\ref{discretizedPDEPOP}) has a unique optimal solution. Therefore, an optimal solution of the SDP (\ref{sparseSDP}) is an approximation of the optimal solution of the POP (\ref{discretizedPDEPOP}).

\subsection{Discrete approximation}

From (\ref{conv}), asymptotic convergence of minimum and minimizer of $\text{SDP}_w$ to minimum and minimizer of POP (\ref{genPOP}) are guaranteed under
uniqueness of the optimal solution and compactness 
of the feasible set of the POP. However, recall that the dimension $\tilde{n}$ of POP (\ref{genPOP}) is $N_x\, N_y$, and so, 
if $w$ is much larger than $w_{\min}$, the resulting $\text{SDP}_w$ in (\ref{sparseSDP}) 
is untractable for general purpose SDP solvers.
Therefore, we have to restrict ourselves to $w\in\{w_{\min},w_{\min}+1 \}$ for medium scale $N$,
and we cannot expect that $\text{SDP}_{w_{\min}}$ or $\text{SDP}_{w_{\min}+1}$ provide accurate approximations to a solution of a nonlinear differential equation.
In order to improve accuracy, we may apply additional locally convergent optimization techniques. 
For instance sequential quadratic programming (SQP) 
can be applied to (\ref{discretizedPDEPOP}), starting from the solution of $\text{SDP}_w$ as initial guess. Since the solution of the SDP relaxation is an approximation for the global optimizer, it is a systematic choice for an initial point for local methods which does not require any a priori information. Combining the sparse SDP relaxation with SQP (or another local method) is summarized in the scheme:

\begin{method} {\it The SDPR method}
\begin{enumerate}
\item 
Choose a discretization $(N_x,N_y)$ for the differential equation problem. 
\item
Choose a relaxation order $w$ and an objective function $F$. Apply $\text{SDP}_w$ to (\ref{discretizedPDEPOP}) and obtain its solution $\tilde{u}$ and the lower bound $\min\left(\text{SDP}_w\right)$ for $\min\left(\text{POP}\right)$.
\item
Apply sequential quadratic programming (SQP) to (\ref{discretizedPDEPOP}) with $\tilde{u}$ as initial guess, and obtain $u$ as discrete approximation to a solution of the differential equation problem.
\end{enumerate}
\label{sdprm}
\end{method}

\section{Maximum entropy estimation}
\label{secMaxEntropy}

In this section we briefly introduce the {\it maximum entropy estimation} of \cite{borwein,lassCDC,lassMPP}, our second tool to find smooth approximations for solutions of nonlinear differential equations. The maximum entropy estimation is concerned with the following problem: Let $u\in L_1(\Omega)$ be nonnegative and partially 
known by the finite vector $m$ of moments up to order $M$ of the associated Borel measure $d\mu:=u\, dx \, dy$ on $\Omega$.
From the only knowledge of $m$, find an estimate $u_M\in L_1(\Omega)$ such that all moments of order up to $M$ of the measure $d\mu_M:=u_M \, dx\, dy$ match those of $d\mu$ and analyze the asymptotic behavior of $u_M$ as $M\rightarrow\infty$.  
An elegant method consists in finding the estimate $u_M$ that maximizes the {\it Boltzmann-Shannon} entropy functional
\[u\mapsto\int_\Omega u\,\ln u\,dxdy,\qquad u\in L_1(\Omega).\]
In this case, the optimal estimate $u_M^\star$ is given by
\begin{equation}
u_M^{\star}(x,y)\equiv \tilde{u}_M^{\star}(v^\star,x,y) = \exp \sum_{\begin{subarray}{c}0\leq i,j\leq M\\ i+j \leq M\end{subarray}} v_{i,j}^\star x^i y^j,
\label{maxEntEstimate}
\end{equation}
where $v^\star\in\RR^{\mid m\mid}$ is a (global) optimizer of the {\it convex} finite-dimensional optimization problem:
\begin{equation}
\max_{v\in\RR^{\mid m\mid}} \langle m,v \rangle - \int_{\Omega} u_M(v,x,y)\, dx\, dy.
\label{concaveOptProb}
\end{equation}
The optimization problem (\ref{concaveOptProb}) can be solved by first or second order methods like Newton's method or SQP. For these methods, gradient and Hessian of the objective function in (\ref{concaveOptProb}) must be computed. In the case the domain $\Omega$ is simple, we may compute them by quadrature and cubature formulas. For more difficult domains one may use the procedure described in \cite{lassCDC,BDL}. Concerning the behavior of 
$u_M^\star$ as $M\to\infty$ and its relationship with $u$, one has the following weak convergence result from \cite{borwein}.
\begin{prop}
\label{convergence}
Let $u_M^\star$ be obtained  from an optimal solution of (\ref{concaveOptProb}).
Then, as $M\to\infty$,
\[\int_\Omega\psi(x,y)\,u_M^\star(x,y)\,dxdy
\:\to\:\int_\Omega\psi(x,y)\,u(x,y)\,dxdy,\]
for every bounded measurable function $\psi:\Omega\to\RR$ which is continuous 
almost everywhere.
\end{prop}
However, the {\it pointwise} convergence $u_M^\star(x,y)\to u(x,y)$ 
does not hold in general.

\section{Smooth approximations to solutions of differential equations}
\label{secSmooth}
In this section we show how to combine the SDPR method and maximum entropy estimation to obtain smooth approximations for solutions of linear and nonlinear differential equations. Our discussion focuses on the 2-dimensional case, but the 1-dimensional case is covered analogously. Let $u\in\D$ be a solution of (\ref{PDEproblem}) and assume without loss of generality that $\text{lbd}\geq0$ so that $u$ is nonnegative on $\Omega$.
For $\text{lbd}<0$ define $\tilde{u}:=u-\text{lbd}$ and apply the outlined procedure to the new function $\tilde{u}$. Associated with $u$, let 
$d\mu(x,y) := u\, dx \, dy$ be the finite Borel measure on $\Omega$
with moment vector $m=(m_{i,j})$ of all moments up to order $M$:
\begin{equation*}
m_{i,j}:=\int_{\Omega} x^i y^j d\mu(x,y)=\int_{\Omega} x^i y^j u(x,y)dx\, dy,
\end{equation*}
for $(i,j)\in\NN^2$ with $i+j\leq M$. For linear PDEs, a hierarchy of tightening lower and upper bounds for the components $m_{i,j}$ of the moment vector can be obtained by solving a sequence of SDPs as proposed in \cite{bertsimasCaramanis}. However, this approach cannot be applied in the case of nonlinear PDEs. Our strategy
is to use a discrete approximation $(u_{k,l})_{1\leq k \leq N_x, 1\leq l \leq N_y}$ of a solution of the PDE to approximate the moments $m_{i,j}$ by:
\begin{equation}
m_{i,j}^\Delta = \sum_{k=1}^{N_x}\sum_{l=1}^{N_y} x_k^i y_l^j u_{k,l}\Delta x \Delta y,
\label{momApprox}
\end{equation}
for $(i,j)\in\NN^2$ with $i+j\leq M$, where $\Delta x:=\frac{x_{\max}-x_{\min}}{N_x-1}$, $\Delta y:=\frac{y_{\max}-y_{\min}}{N_y-1}$. If the discretization $(N_x, N_y)$ is sufficiently fine and $(u_{k,l})_{k,l}$ is a close approximation of $u$, then we expect $m_{i,j}^\Delta$ to be a good approximation of $m_{i,j}$ for all $(i,j)\in\NN^2$. Thus, we can apply maximum entropy estimation to the vector $m^\Delta$ to obtain a smooth approximation of $u$. This idea is formalized in the following algorithm for obtaining smooth approximations of solutions of linear or nonlinear PDE problems.

\begin{method} {\it The smooth SDP approximation method}\\
Given a PDE problem of form (\ref{PDEproblem}).
\begin{enumerate}
\item
Choose a discretization $(N_x,N_y)$, relaxation order $w$ and objective $F$, and apply the SDPR method to obtain a discrete approximation $(u_{k,l})_{k,l}$ to a solution of the PDE problem. If not given in the formulation of the PDE problem, impose lower and upper bounds, lbd and ubd for $u$.
\item
Choose a moment bound $M\in\NN$ and use $(u_{k,l})_{k,l}$ to calculate $m^\Delta$ as 
in (\ref{momApprox}).
\item
Apply maximum entropy estimation to $m^\Delta$ and obtain vector $v^\star\in\RR^{\mid m^\Delta\mid}$, optimal solution of (\ref{concaveOptProb}).
\item
Obtain the approximation $u_M^{\star}$ with $u_M^{\star}(x,y)=\sum_{i,j}\exp(v_{i,j}^\star x^iy^j)$ for a solution $u$ of the PDE problem (\ref{PDEproblem}) on $\Omega$.
\end{enumerate}
\label{sSDPam}
\end{method}
As $M\to\infty$, $u_M^{\star}\rightarrow u$ weakly (see Proposition \ref{convergence})
but not pointwise, i.e., one cannot guarantee $u_M^{\star}(x,y)\rightarrow u(x,y)$ 
on $\Omega$. Nevertheless, as reported in \cite{lassCDC} the maximum entropy estimation may provide accurate pointwise approximation of the unknown function to be recovered on certain segments of the domain $\Omega$. We next illustrate on a variety of PDE problems and OCPs, that indeed good pointwise approximation can be obtained on 
some parts of the domain $\Omega$.

\section{Numerical Experiments}
\label{secNumExp}
We illustrate the potential of Method \ref{sSDPam} on a variety of ODE and PDE problems. As an implementation of the sparse SDP relaxations we apply the software SparsePOP of \cite{wkkm2} and as an implementation of sequential quadratic programming (SQP) in the SDPR method and in solving the optimization problem (\ref{concaveOptProb}) we apply the Matlab Optimization Toolbox commands {\tt fmincon} and {\tt fminunc}, respectively. As we restrict ourselves to ODEs and PDEs with rectangular domains $\Omega$ we can apply standard quadrature and cubature formulas to compute the gradient and Hessian for (\ref{concaveOptProb}). Thus, we apply the Matlab commands {\tt trapz} in the one-dimensional case and {\tt dblquad} in the two-dimensional case, respectively.
In order to evaluate the quality of the smooth approximation provided by Method \ref{sSDPam} we define the {\it average error} $\bar{\epsilon}_u(M):= \Delta x\sum_i u(x_i)-u_M^\star(x_i) $ and the {\it maximum pointwise error} $\epsilon_u^{\max}(M):=\max_i \mid u(x_i)-u_M^\star(x_i)\mid$.

\subsection{Linear differential equations}
As first test problems for our approach we consider a linear ODE and a linear PDE from \cite{bertsimasCaramanis} and compare (a) the discretized moment approximations obtained by Method \ref{sSDPam} to the bounds obtained in \cite{bertsimasCaramanis}, and (b) the smooth approximation $u_M^\star$ to the known analytic solution. 

\subsubsection{Linear ODE}

The linear ODE is given by,
\begin{equation}
\begin{array}{ll}
u''(x) + 3u'(x) + 2u(x) = 0 & \forall \, x\in [0,1],\\
u'(0) = -2e^2, \,\, u'(1) = -2.
\end{array}
\label{linearODE} 
\end{equation}
For this problem the unique solution is $u(x)=e^2e^{-2x}$. We apply Method \ref{sdprm} with $F=-\sum_i u_i$, $N=2000$ and $w=1$. We calculate the approximate moments $m^\Delta$ for $M=40$ and compare them to approximate moments for $m$ derived in \cite{bertsimasCaramanis} by contracting lower and upper bounds in Table \ref{resLinearODE}. All moments coincide up to the fourth digit.

\begin{table}
\begin{center}
\begin{tabular}{|r|l|l|}
\hline $i$ & $m_i^{\Delta}$ & $m_i^{BC}$ \\ 
\hline 0 & 3.1942 & 3.1945 \\
\hline 1 & 1.0957 & 1.0973 \\
\hline 10 &  0.1086 & 0.1088 \\
\hline 20 & 0.0524 & 0.0524 \\
\hline 30 & 0.0345 & 0.0344 \\
\hline 40 & 0.0258 & 0.0256 \\
\hline
\end{tabular}
\caption{Approx. moments $m^\Delta$ for $N=2000$ compared to approx.  $m^{BC}$ derived from lower and upper bounds in \cite{bertsimasCaramanis} for linear ODE (\ref{linearODE})}
\label{resLinearODE}
\end{center}
\end{table}

We apply Method \ref{sSDPam} for $M\in\{1,\ldots,5\}$ and report the resulting vector $v_M^\star$ in Table \ref{resLinearODEvstar}.

Note that the actual solution of the ODE problem corresponds to a parameter vector $v_{\text{opt}}=(2, -2, 0,\ldots)$. As reported in Table \ref{resLinearODEvstar}, $(v_0^\star(M), v_1^\star(M))\approx (2,-2)$, and even though $v_i^\star(M) \neq 0$ for $i>1$, the maximum pointwise error is quite small. However, $(v_0^{\star}(M),v_1^{\star}(M))$ does not converge to $(2,-2)$ since pointwise convergence is not guaranteed. In fact the higher moment terms $v_i^{\star}(M)$ for $i>1$ counterbalancce the difference between  $(v_0^{\star}(M),v_1^{\star}(M))$ and $(2,-2)$.

\begin{table}
\begin{center}
\begin{tabular}{|l|r|r|r|r|r|}
\hline $M$ & 1 & 2 & 3 & 4 & 5\\
\hline $v_{0}^\star$ & 2.0043 & 2.0042 & 2.0040 & 2.0032 & 2.0039\\
$v_1^\star$ & -2.0135 & -2.0083 & -1.9922 & -1.9744 & -1.9528\\
$v_2^\star$ & & -0.0083 & -0.0938 & -0.1998 & -0.3205 \\
$v_3^\star$ & & & 0.0851 & 0.2886 & 0.3269 \\
$v_4^\star$ & & & & -0.1180 & 0.2018 \\
$v_5^\star$ & & & & & -0.2694 \\
\hline $\bar{\epsilon}_u(M)$ & 0.0042 & 0.0044 & 0.0042 & 0.0039 & 0.0044 \\
\hline $\epsilon_u^{\max}(M)$ & 0.0242 & 0.0228 & 0.0228 & 0.0222 & 0.0326 \\
\hline
\end{tabular}
\caption{Smooth SDP approx. for linear ODE (\ref{linearODE})\label{resLinearODEvstar}}
\end{center}
\end{table}

\subsubsection{Linear PDE}

The linear PDE is given by
\begin{equation}
u_{xx}(x,y) + u_{yy}(x,y) - 3e^{x+y} = 0 \:\: \forall \, (x,y)\in [0,1]^2,
\label{linearPDE} 
\end{equation}
where the boundary conditions are set up such that $u(x,y)=e^{x+y}$ is the unique solution of the PDE problem. We apply the SDPR method with $F(u)=-\sum_{k,l}u_{k,l}$, $w=1$ and $N_x=N_y=100$, compute $m^\Delta$ for all moments of order up to 3 and compare these moments to the exact moments derived from the known solution and to the best lower and upper bounds for the moments from Table 7 in \cite{bertsimasCaramanis}. See Table \ref{resLinearPDE} for the results. Our approach provides some approximations for the moments that  is well within the bounds LB and UB from \cite{bertsimasCaramanis}.

\begin{table}
\begin{center}
\begin{tabular}{|l|l|l|l|l|}
\hline $(i,j)$ & $m_{i,j}^\Delta$ & $m_{i,j}$ & LB & UB \\
\hline (0,0) & 2.9512 & 2.9525 & 2.9235 & 3.1707 \\
\hline (1,0) & 1.7287 & 1.7183 & 1.6944 & 1.7742 \\
\hline (1,1) & 1.0121 & 1.0000 & 0.9847 & 1.0130 \\
\hline (2,0) & 1.2504 & 1.2342 & 1.1123 & 1.3088 \\
\hline (2,1) & 0.7319 & 0.7183 & 0.7151 & 0.7458 \\
\hline (2,2) & 0.5292 & 0.5159 & 0.4948 & 0.5456 \\
\hline (3,0) & 0.9869 & 0.9681 & 0.8244 & 1.0478 \\
\hline (3,1) & 0.5775 & 0.5634 & 0.5054 & 0.5818 \\
\hline (3,2) & 0.4176 & 0.4047 & 0.3874 & 0.5818 \\
\hline (3,3) & 0.3295 & 0.3175 & 0.3017 & 0.3399 \\
\hline
\end{tabular}
\caption{Approx. moments for $N_x=N_y=100$ compared to exact moments and bounds from \cite{bertsimasCaramanis}
for linear PDE (\ref{linearPDE})\label{resLinearPDE}}
\end{center}
\end{table}

Also for this problem, there is a vector $v_{\text{opt}}=(0,1,1,0,\ldots)$ corresponding to the solution $u$. When applying Method \ref{sSDPam} with $M\in\{2,3\}$ we obtain $v^\star(M)$ reported in Table \ref{resLinearPDEvstar} with errors $\bar{\epsilon}_u(2)=0.0677$ , $\bar{\epsilon}_u(3)=0.0682$, $\epsilon_u^{\max}(2)=1.0276$ and $\epsilon_u^{\max}(3)=1.3267$.

\begin{table}
\begin{center}
\begin{tabular}{|l|rrrrrr|}
\hline $(i,j)$ & (0,0) & (1,0) & (0,1) & (2,0) & (1,1) & (0,2)\\
\hline $v_{i,j}^\star(2)$ & 0.064 & 0.709 & 0.710 & 0.329 & 0.099 & 0.329 \\
\hline $v_{i,j}^\star(3)$ & 0.013 & 0.959 & 0.959 & -0.213 & 0.160 & -0.213 \\
\hline ${v_{\text{opt}}}_{i,j}$ & 0 & 1 & 1 & 0 & 0 & 0 \\
\hline
\end{tabular}
\caption{Smooth SDP approx. for linear PDE (\ref{linearPDE})\label{resLinearPDEvstar}}
\end{center}
\end{table}
For both linear differential equations we obtain accurate approximations of the moments for the measures associated with the unique solutions. Moreover, unlike \cite{bertsimasCaramanis}, we can exploit these moment approximation to find smooth approximations for the actual solutions.

\subsection{A nonlinear elliptic PDE}

\begin{figure*}[thpb]
\centering
\includegraphics[width=0.3\textwidth, height=0.17\textheight]{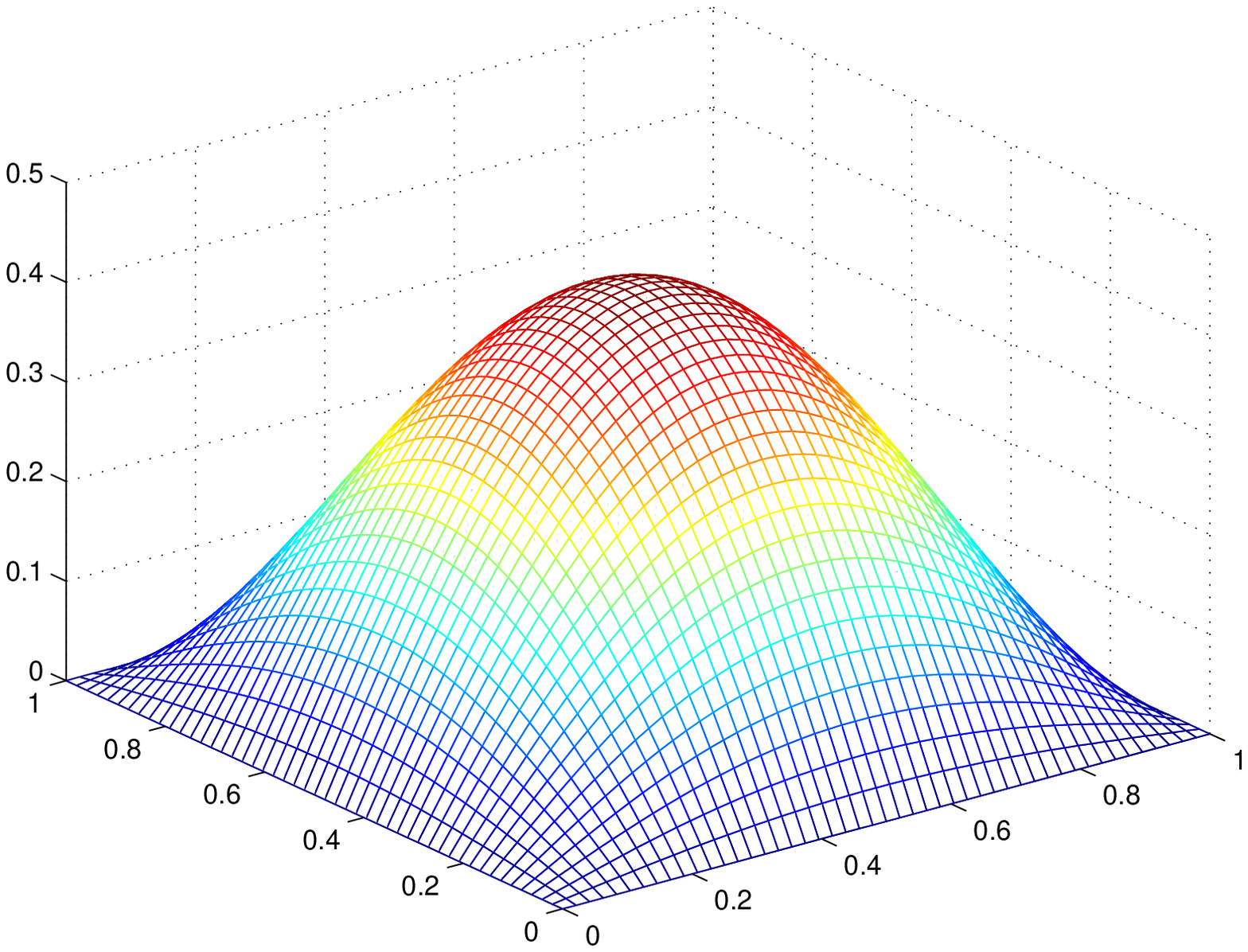} 
\includegraphics[width=0.3\textwidth, height=0.17\textheight]{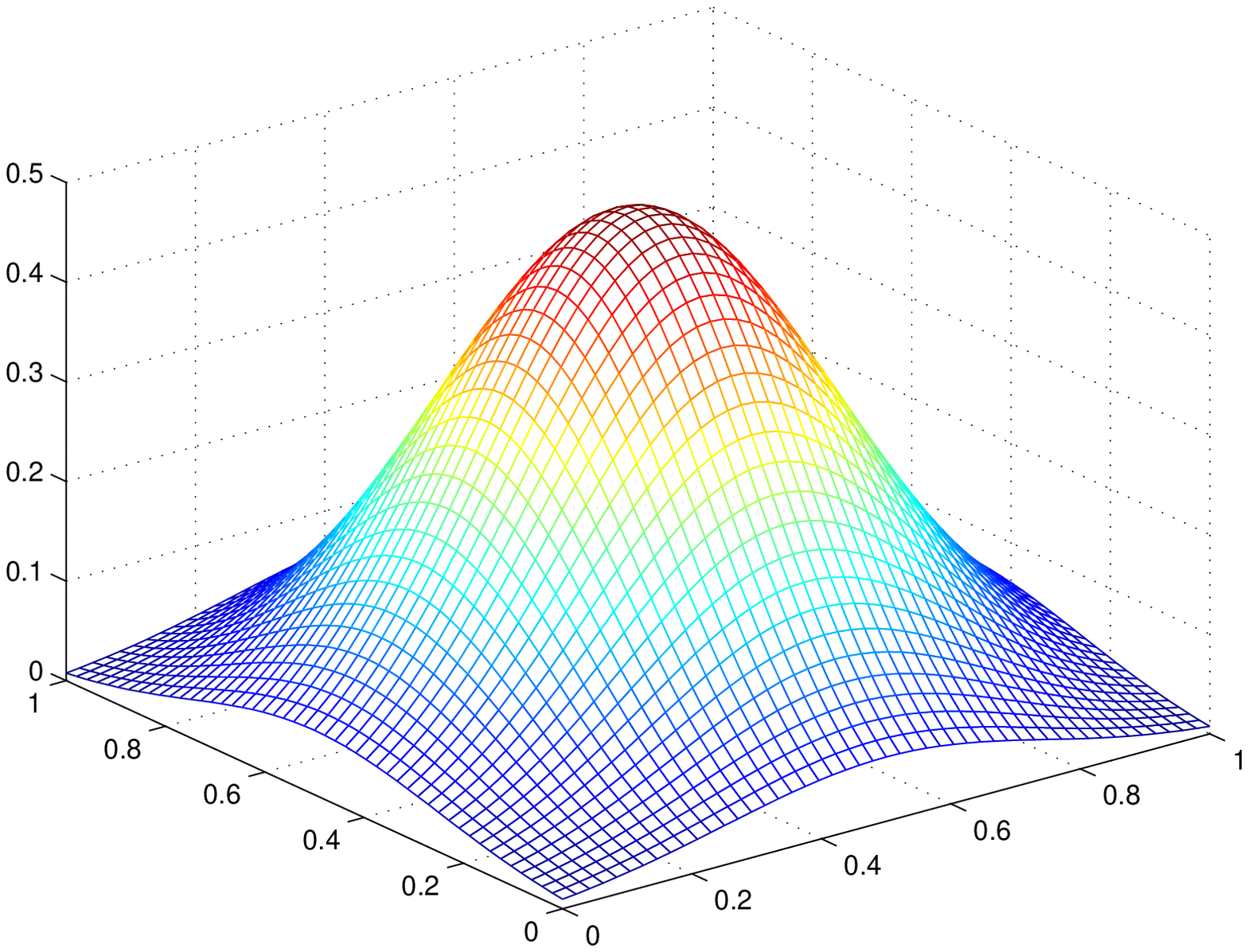} 
\includegraphics[width=0.3\textwidth, height=0.17\textheight]{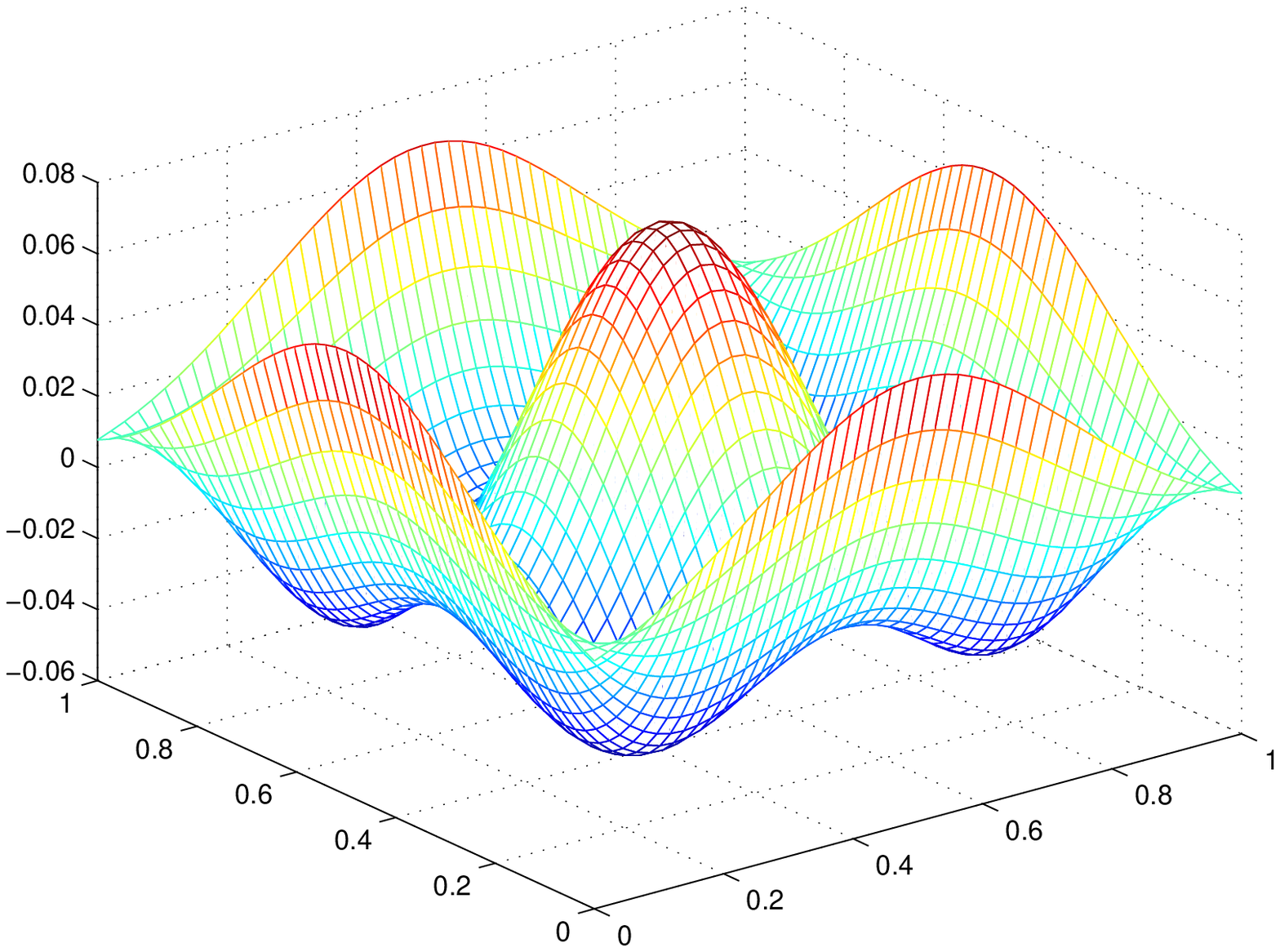} 
\caption{Discrete approx. (left), smooth approx. (center) for the same grid, and pointwise difference between the two (right, notice the vertical axis scale) for nonlinear PDE (\ref{pdeBifur})\label{bifurFig}}
\end{figure*}

As a first nonlinear problem, consider the elliptic PDE:
\begin{equation}
\begin{array}{ll}
u_{xx} + u_{yy} + 22u(1-u^2) =0 & \text{on } [0,1]^2,\\
u = 0& \text{on } \partial [0,1]^2,\\
0\leq u\leq 1 & \text{on } [0,1]^2.
\end{array}
\label{pdeBifur}
\end{equation}
This problem is well known to have a nontrivial positive solution. We apply Method \ref{sSDPam} with $F(u):=-\sum_{k,l}u_{k,l}$, $w=2$, $N_x=N_y=49$ and $M\in \{2,3\}$, and obtain an accurate discrete approximation $(u_{i,j})_{i,j}$ and two smooth approximations $u_2^\star$ and $u_3^\star$ with errors $(\epsilon_u^{\max}(2),\bar{\epsilon}_u(2)) = (0.081, -0.002)$ and $(\epsilon_u^{\max}(3),\bar{\epsilon}_u(3)) = (0.076, -0.002)$, respectively. In Figure \ref{bifurFig} it is illustrated how accurate the smooth approximation $u_3^\star$ resembles the shape of the discrete approximation $(u_{i,j})_{i,j}$.

\subsection{Reaction diffusion equation} 

\begin{figure*}[thpb]
\centering
\includegraphics[width=0.45\textwidth]{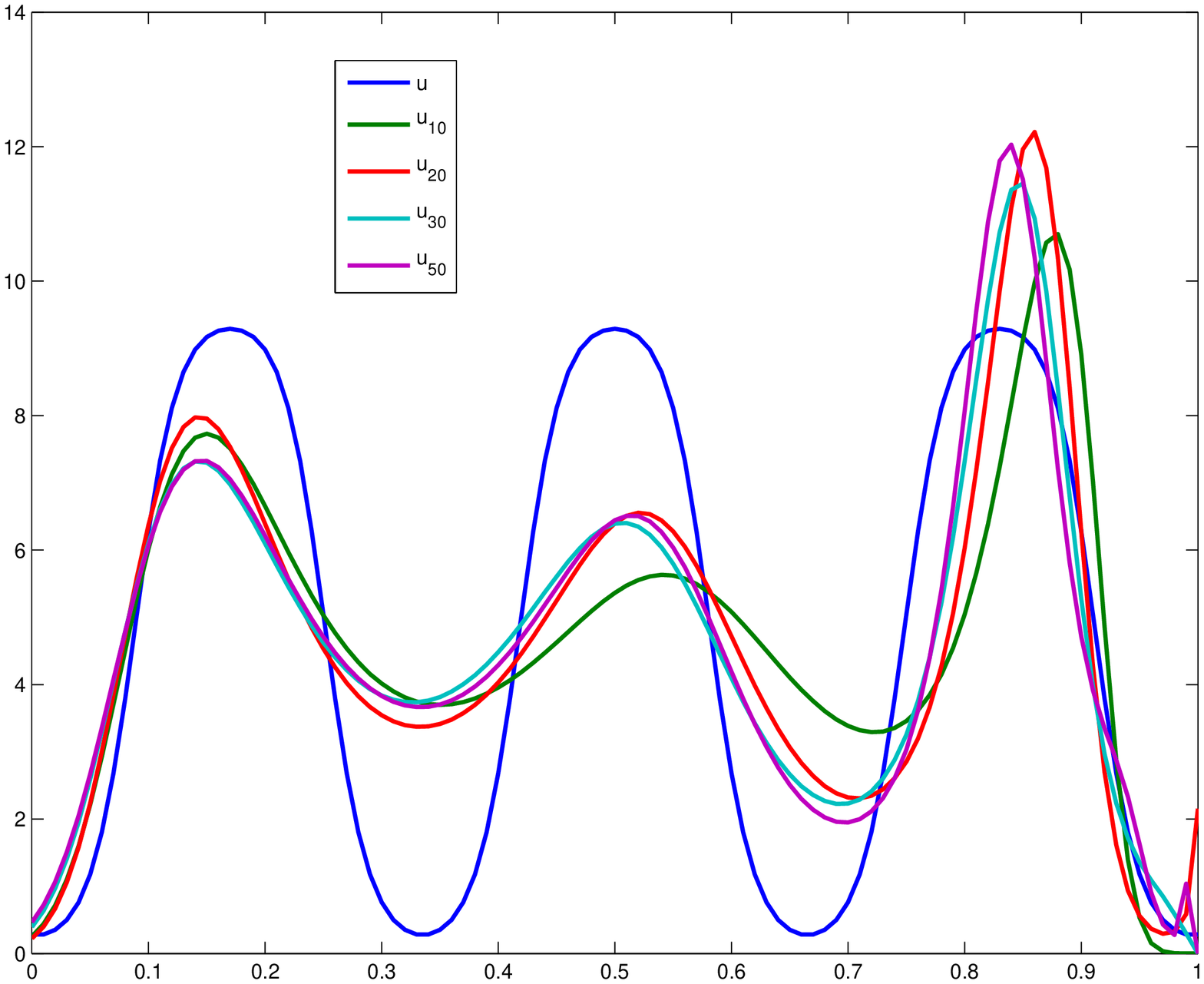} 
\includegraphics[width=0.45\textwidth]{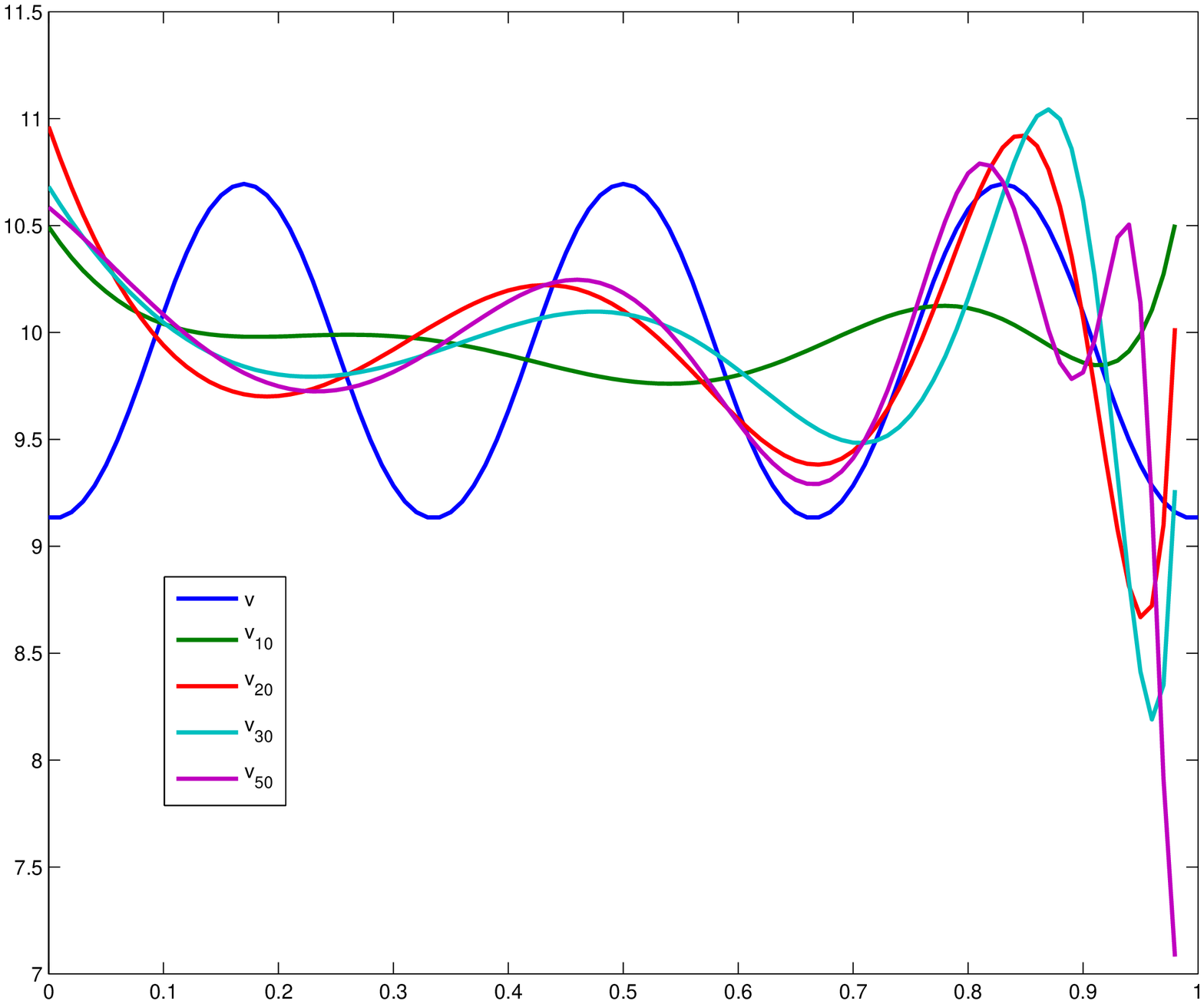} 
\caption{Smooth approx. for $u$ (left) and $v$ (right) for $M\in\{10,20,30,50\}$ for nonlinear ODE (\ref{mimuraODE})
\label{mimuraFig}}
\end{figure*}

A challenging ODE problem which is known to have many solutions is given as in \cite{mimura} by :
\begin{equation}
\begin{array}{l@{\:}l}
\frac{1}{20} \; {u''}+ \frac{1}{9} \left( 35 + 16u - u^2 \right)\; u - u\, v = 0 & \text{on } [0,5], \\
4 {v''} - \left( 1 + \frac{2}{5}v\right)\; v + u\, v = 0 & \text{on } [0,5],\\
u'(0) = u'(5) = v'(0) = v'(5) = 0,\\
0\leq u,v  \leq 14 & \text{on }[0,5]. 
\end{array}
\label{mimuraODE}
\end{equation}
To ensure numerical stability we scale the domain $[0,5]$ of (\ref{mimuraODE}) to $[0,1]$ before applying Method \ref{sSDPam}, as $\exp(x^i)$ gets very large for $|x|>1$. Problem (\ref{mimuraODE}) involves two functions $u$ and $v$, thus we need to apply Method \ref{sSDPam} twice, once for $u$ and once for $v$, in order to obtain smooth approximations for both functions. Note that the SDPR method only needs to be applied once and provides a discrete approximation $(u_{i},v_i)_i$, maximum entropy estimation needs to be applied twice. We apply Method \ref{sSDPam} for $F(u,v)=-u_{\lceil\frac{N}{2}\rceil}$ , $N=100$, $w=3$ and $M\in\{ 10,20,30,50 \}$. Solutions for $u$ and $v$ and the derived smooth approximations are pictured in Figure \ref{mimuraFig} and the errors are reported in Table \ref{mimuraError}. Even though the average and the maximum pointwise approximation errors are not small, even for large $M$, we observe that the pointwise approximation gets more and more accurate on certain segments of the domain - in this case the third part of the interior, where the third peak of $u$ and $v$ occurs. The pointwise approximation is not good - or does even get worse - near the boundary of $[0,1]$.

\begin{table}
\begin{center}
\begin{tabular}{|l|r|r|r|r|}
\hline $M$ & $\epsilon_u^{\max}(M)$ & $\bar{\epsilon}_u(M)$ & $\epsilon_v^{\max}(M)$ & $\bar{\epsilon}_v(M)$ \\
\hline 10 & 4.11 & 0.0011 & 2.10 & 0.1080\\
20 & 3.84 & 0.0143 & 6.66 & 0.1371 \\
30 & 3.48 & 0.0012 & 9.61 & 0.1547 \\ 
50 & 3.39 & 0.0085 & 23.69 & 0.2658\\
\hline
\end{tabular}
\caption{Errors for nonlinear ODE (\ref{mimuraODE})
\label{mimuraError}}
\end{center}
\end{table}

\subsection{Control of production and consumption}

\begin{figure*}[thpb]
\centering
\includegraphics[width=0.45\textwidth]{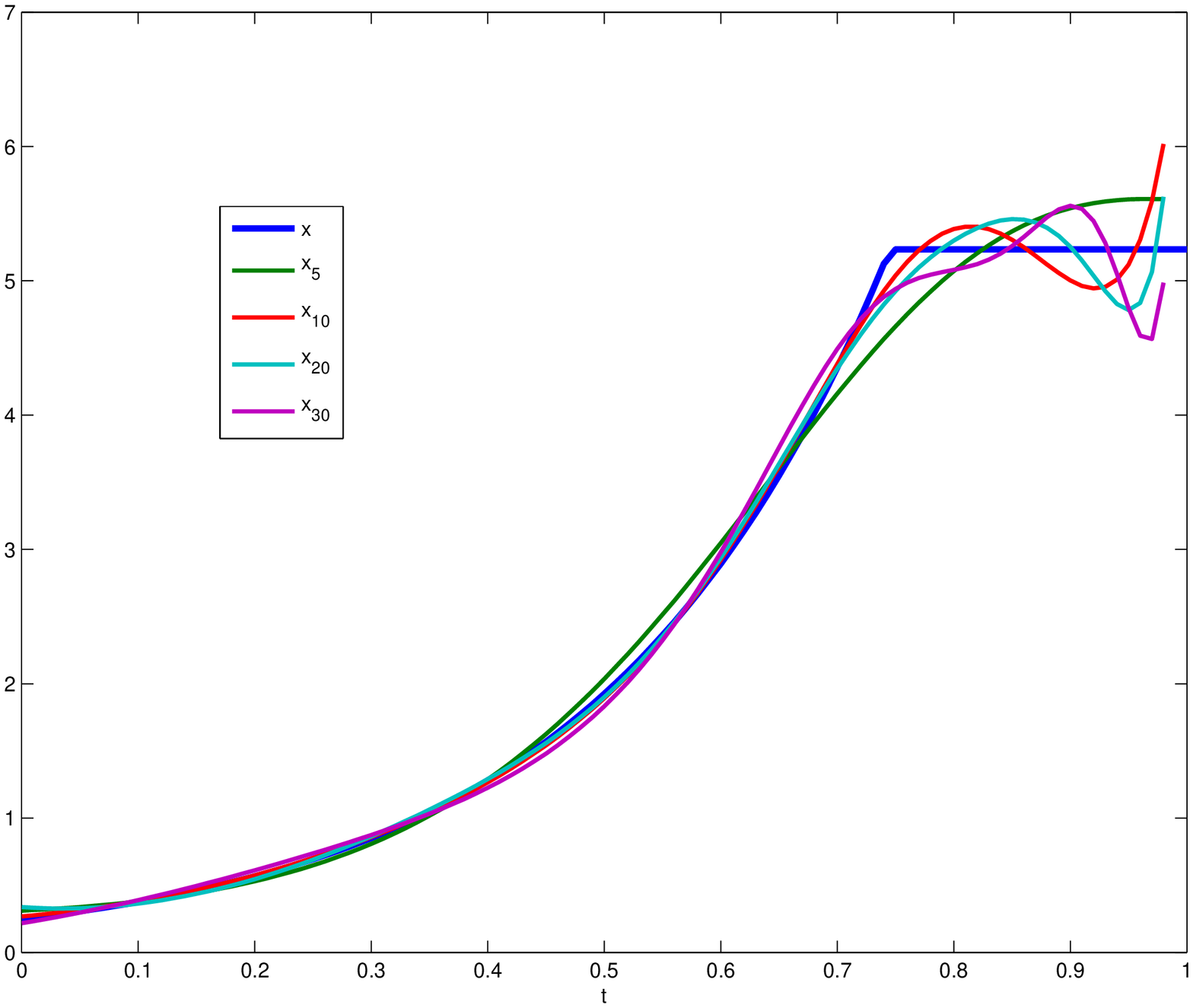} 
\includegraphics[width=0.45\textwidth]{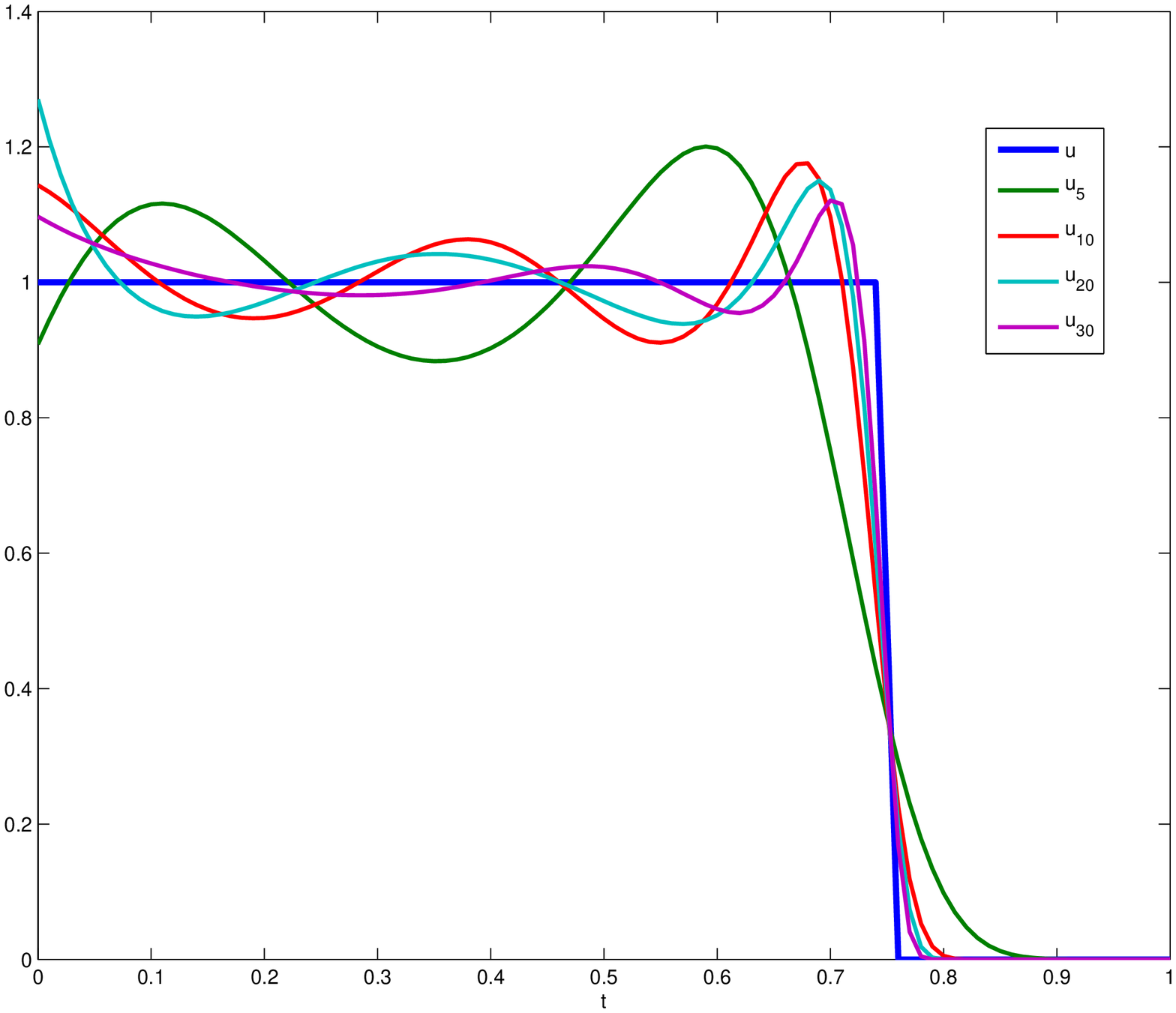} 
\caption{Smooth approx. for $x$ (left) and $u$ (right) for $M\in\{5,10,20,30\}$ for OCP (\ref{controlProd}).}
\label{prodFig}
\end{figure*}

A first example of an optimal control problem is given by
\begin{equation}
\begin{array}{lll}
\min & -\int_0^T (1-u(t))x(t) \, dt \\
\text{s.t.} & \dot{x}(t)  = u(t) x(t) & \forall t \in \left[ 0,T\right] , \\
& x(0) = x_0, \\
& 0\leq x(t) \leq 1 & \forall t \in \left[ 0,1 \right], \\
& 0\leq x(t) \leq 10 & \forall t \in \left[0,10 \right],
\end{array}
\label{controlProd}
\end{equation}
where $T>1$ fixed. For this simple problem the optimal control law is given by 
\begin{equation*}
u^{\star}(t) =\begin{cases} 1 & \text{if } 0\leq t \leq T-1,\\ 0 & \text{if } T-1< t \leq T. \end{cases}
\end{equation*}
We choose $x_0=0.25$, $T=4$ and $N=100$, and apply Method \ref{sSDPam} after scaling the domain to $[0,1]$ for $M\in\{ 5,10,20,30 \}$. In this case the objective function $F$ is given by a discretization of the objective function of (\ref{controlProd}). As in the case of ordinary and partial differential equations, we observe that a fairly good pointwise approximations of both, optimal control law and corresponding trajectory, are obtained on the interior of the domain. See Figure \ref{prodFig}.

\subsection{The double integrator}

\begin{figure*}[thpb]
\centering
\includegraphics[width=0.3\textwidth]{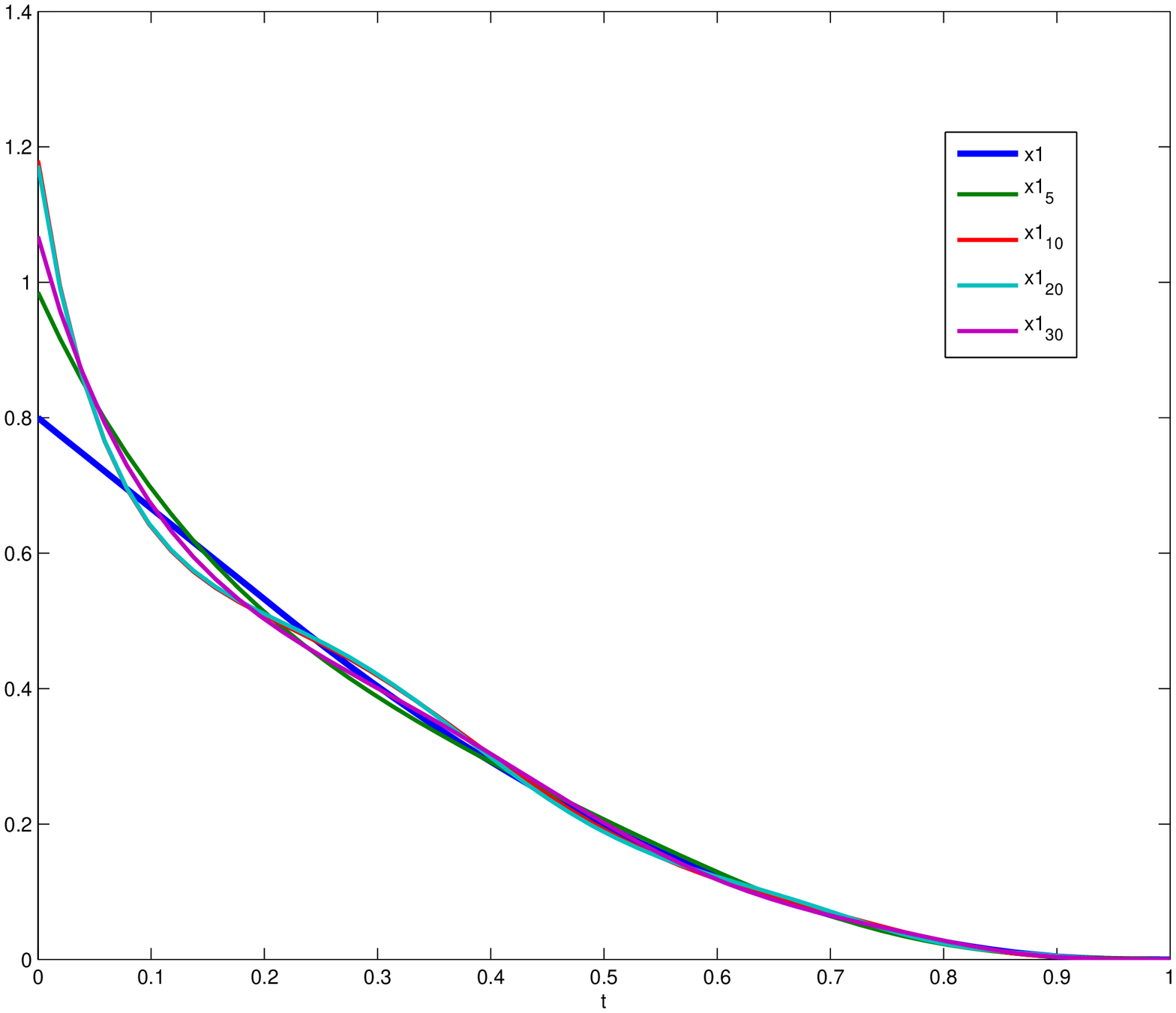} 
\includegraphics[width=0.3\textwidth]{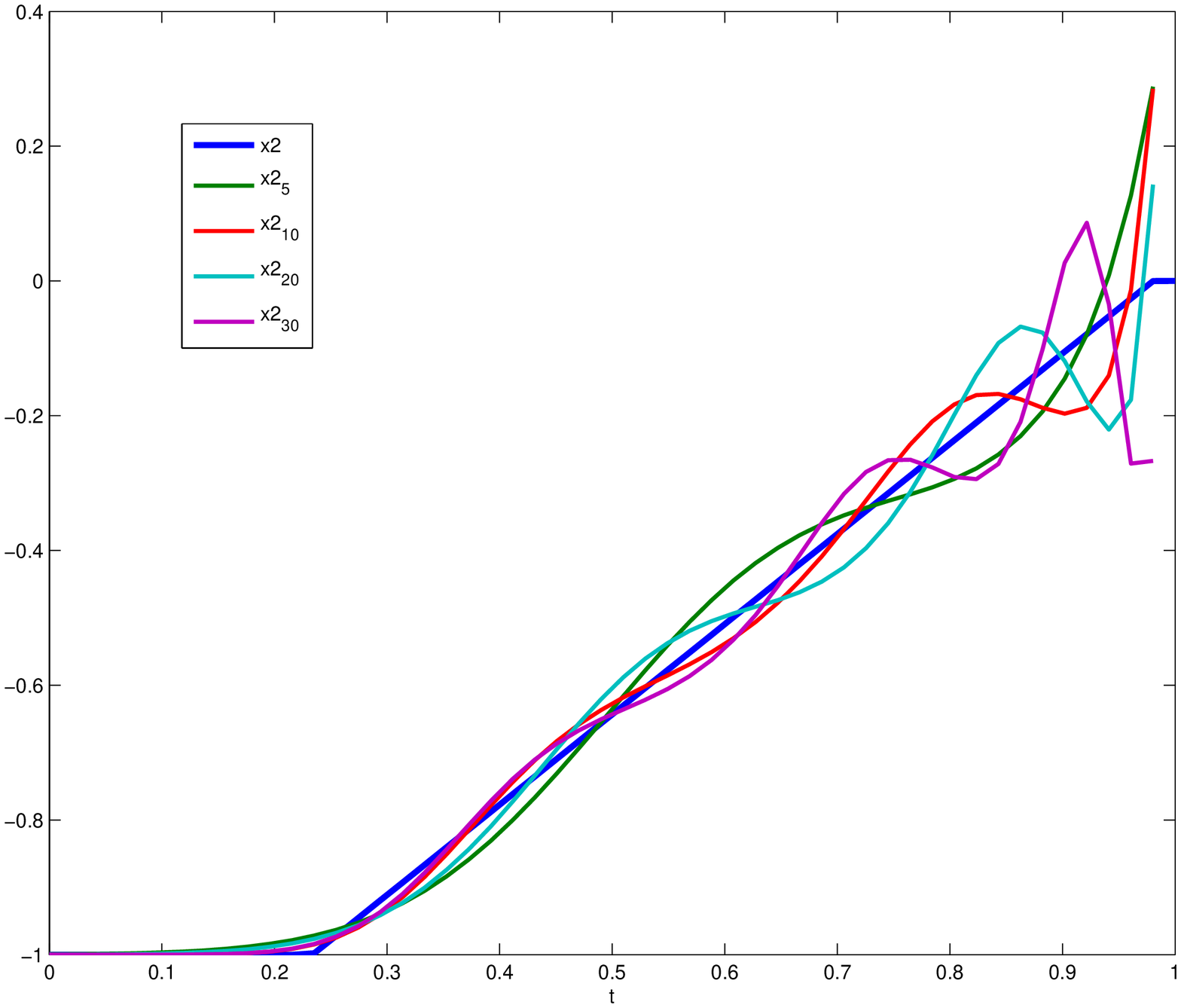} 
\includegraphics[width=0.3\textwidth]{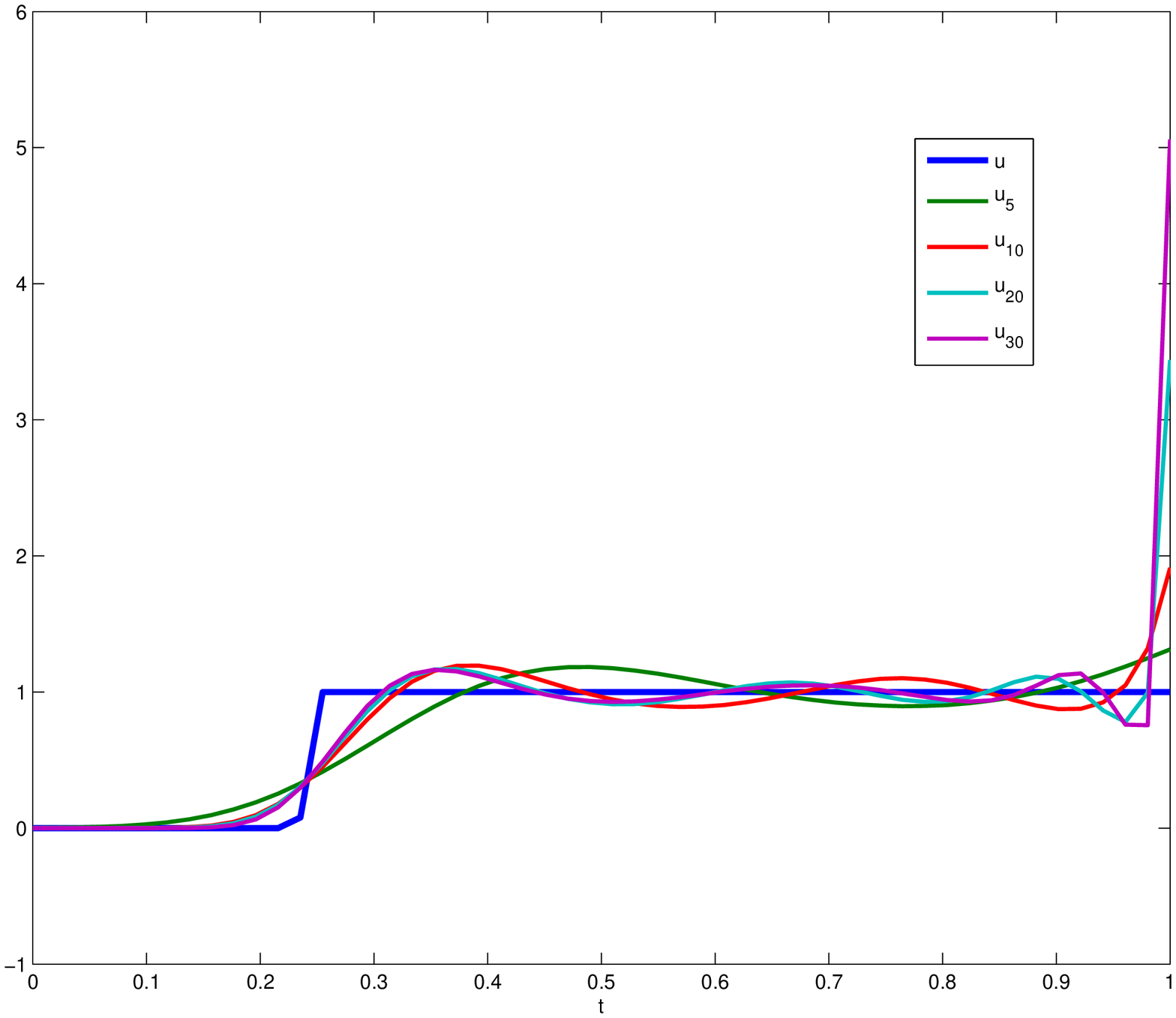} 
\caption{Smooth approx. for $x_1$ (left), $x_2$ (center) and $u$ (right) for $M\in\{5,10,20,30\}$ forOCP (\ref{controlDouble}).}
\label{doubleFig}
\end{figure*}

Another interesting control problem, which has been discussed in \cite{lhpt}, is given by
\begin{equation}
\begin{array}{lll}
\min & T \\
\text{s.t.} &\dot{x}_1(t) = x_2(t) & \forall \, t\in \left[ 0,T \right], \\
& \dot{x}_2(t) = u(t) & \forall \, t \in \left[ 0,T\right],\\
& x(0) = x_0\in\RR^2,\\
& x(T) = (0,0),\\
& -1\leq u(t) \leq 1 & \forall \, t \in [0,T],\\
& -1\leq x_1(t), x_2(t) \leq 10 & \forall \, t \in [0,T].
\end{array}
\label{controlDouble}
\end{equation}
After scaling the domain to $[0,1]$, we apply Method \ref{sSDPam} with $N=50$, $\omega=3$, $x_0=(0.8,-1)$ and $M\in\{5,10,20,30\}$. Since $T$ is not fixed, it is treated as an additional variable, i.e., we consider the polynomial optimization problem with variable $({x_1}_1,\ldots, {x_1}_N,{x_2}_1,\ldots,{x_2}_N,u_1,\ldots,u_N,T)$. We observe that the pointwise approximation of optimal control and trajectories on the interior of the domain gets better and better for increasing moment order, see Figure \ref{doubleFig}. Note, unlike the moment based method in \cite{lhpt} which yields bounds for the optimal value of the control problem, we obtain discrete and smooth approximations to the optimal control and the corresponding trajectories by Method \ref{sdprm} and Method \ref{sSDPam}, respectively.

\section{Concluding Remarks}

We introduced a novel technique to derive smooth approximations for solutions of systems of differential equations and optimal control problems, which is based on sparse SDP relaxations and the maximum entropy estimation method. As demonstrated on some examples of nonlinear ordinary differential equations, partial differential equations and optimal control problems, this technique is promising to obtain accurate pointwise approximations of a solution of a differential equation on segments of its domain. It would be interesting to characterize regions of the domain of nonlinear differential equations where accurate pointwise approximations can be guaranteed. Another question is whether a different, less restrictive parametrization than (\ref{maxEntEstimate}) allows better pointwise approximations on larger segments of the domain. 
Of course, other choices of entropy would lead to different results, and it remains a topic of future investigation to analyze, whether or not some choice should be preferred. At this point solving the sparse SDP relaxation is the major computational bottleneck in Method \ref{sSDPam}. Every improvement for approximation accuracy and efficiency of SDP relaxation techniques for polynomial programs and every improvement for SDP solvers will extend the range of the proposed method.

In order to validate these techniques in the present context of PDEs, a comparison with state-of-the-art numerical methods for solving nonlinear PDEs remains to be done. In particular, the maximum entropy technique should be compared with standard and spline
interpolation methods which both provide differentiable approximations that coincide with the discrete solution at each grid point. But, unlike our method, in polynomial interpolation the {\it degree} of the polynomial is directly related to the number of grid points. Also, cubic splines are twice differentiable at the grid points, whereas our approximation is smooth. 
In contrast, and instead of searching for a smooth solution that matches the unknown $u$ at the grid points, one searches for a smooth solution that matches finitely many {\it moments} of the associated measure $d\mu=udx$. And so, an interesting feature is to
obtain a smooth approximation from a limited number of {\it moments} and not from the discrete solution explicitly (still, before using moments, a sufficiently fine discretization is necessary to ensure $m^\Delta\approx m$).
This method is an attempt to deal with the {\it curse of dimensionality} in the numerical analysis of PDEs. Furthermore, we emphasize that the proposed method can be applied to a wide range of problems that encompasses nonlinear partial differential equations and nonlinear optimal control problems. A detailed comparison is however out of the scope of this paper, whose main objective was to pave the way for the development of SDP techniques to solve optimal control design problems for dynamical systems described by polynomial ODEs and PDEs, in the spirit of \cite{lhpt}.



\section*{Acknowledgements}

The research of M. Mevissen was supported by the Doctoral Scholarship of the German Academic Exchange Service. The research of D. Henrion was partly supported by project No.~103/10/0628 of the Grant Agency of the Czech Republic.

\end{document}